\documentclass{amsart}
\usepackage{amsfonts,amssymb,amsmath,amsthm}
\usepackage{url}
\usepackage{enumerate}

\newtheorem{theorem}{Theorem}[section]
\newtheorem{lemma}[theorem]{Lemma}
\newtheorem{corollary}[theorem]{Corollary}
\newtheorem{proposition}[theorem]{Proposition}
\theoremstyle{definition}

\numberwithin{equation}{section}

\def\la{\lambda}

\def\R{\mathbb{R}}

\def\lo{{\Lambda^p_u(w)}}
\def\loin{{\Lambda^{p,\infty}_u(w)}}
\def\llo{{\Lambda^p(w)}}
\def\lloin{{\Lambda^{p,\infty}(w)}}

\def\bdem{\begin{proof}}
\def\edem{\end{proof}}

\def\ds{\displaystyle}

\begin{document}

\author {Elona Agora, Jorge Antezana,  Mar\'{\i}a J. Carro} 

\address{E. Agora, Instituto Argentino de Matem\'atica ``Alberto P. Calder\'on'', 
1083 Buenos Aires, Argentina.}
\email{elona.agora@gmail.com}

\address{J. Antezana, Department of Mathematics, 
Faculty of Exact Sciences, 
National University of La Plata, 
1900 La Plata, Argentina,  \and
Instituto Argentino de Matem\'atica ``Alberto P. Calder\'on'', 
1083 Buenos Aires, Argentina.}
\email{antezana@mate.unlp.edu.ar}

\address{M. J. Carro, Department of Applied Mathematics and Analysis, University of Barcelona, 
 08007 Barcelona, Spain.}
\email{carro@ub.edu}

\subjclass[2010]{26D10, 42A50}
\keywords{Weighted Lorentz spaces, Hardy-Littlewood maximal operator}
\thanks{This work has been partially supported by Grants MTM2010-14946, MTM2013-40985-P and  2014SGR289, and
the University of La Plata Grant 11X681, CO\-NI\-CET Grant PIP-2009-435 and the PICT 2011-436.}

\title{Weak-type boundedness of Hardy-Littlewood maximal operator on weighted Lorentz spaces} 

\begin{abstract} 
The main goal of this paper is to provide a complete characterization of the weak-type boundedness of the Hardy-Littlewood 
maximal operator, $M$, on weighted Lorentz spaces $\Lambda^p_u(w)$, whenever $p>1$. This solves a problem left open 
in \cite{crs:crs}. Moreover, with this result, we complete the program of unifying the 
study of the boundedness of $M$ on weighted Lebesgue spaces  and classical Lorentz spaces, which was initiated in the 
aforementioned monograph.
\end{abstract}

\maketitle

\markboth{Weak-Type Boundedness of the Hardy–Littlewood Maximal Operator on Weighted Lorentz Spaces}
{Elona Agora, Jorge Antezana, and  Mar\'{\i}a J. Carro}

\section{Introduction}

\noindent The classical Hardy-Littlewood maximal operator $M$, is defined by
$$
Mf(x)=\sup_{x\in Q}{\frac{1}{|Q|}}\int_{Q} |f(y)|dy,
$$
where the supremum is taken over all cubes $Q$ containing $x \in \R^d$. This operator is related with several problems in analysis, 
and in some sense it controls the boundedness of many other operators. For these reasons, it has been widely studied in different settings. 

In 1972, Muckenhoupt \cite{m:m} gave the complete characterization of the boundedness of  $M$ on
weighted Lebesgue spaces $L^p(u)$, defined by the set of all Lebesgue measurable
functions $f$ such that
$$
||f||_{L^p(u)}:=\left( \int_{\R^d} |f(x)|^p u(x)dx \right)^{1/p}<\infty,
$$
where $u$ is a positive and locally integrable function on $\mathbb R^d$  (we call it weight).
For $p>1$, the characterization was given in terms of the so called  $A_p$ class of weights \cite{m:m}; that is 
\begin{equation*}\label{defAp}
\sup_{Q}\left(\frac{1}{|Q|}\int_Q u(x)dx\right) \left(\frac{1}{|Q|} \int_Q u^{-1/(p-1)}(x) dx \right)^{p-1}<\infty,
\end{equation*}
where the supremum is considered over all cubes $Q$ of $\R^d$. It was also proved in \cite{m:m} that, if $p > 1$, 
$$
M : L^p(u) \to L^p(u)  \iff M : L^p(u)\to L^{p,\infty}(u)  \iff  u\in A_p,
$$
where the weak-type space $L^{p,\infty}(u)$ is defined through the quasi norm
$$
||f||_{L^{p,\infty}(u)}:= \sup_{t>0} t \, u(\{x\in \R^d: \, |f(x)|>t\})^{\frac 1p}<\infty,
$$
and $\ds u(E)= \int_E u(x)dx$,  for any measurable set $E\subset \R^d$.  
If $p=1$ the only case that makes sense is the weak-type boundedness 
$$
M : L^1(u)\to L^{1,\infty}(u),
$$
characterized by the $A_1$ class of weights defined by 
$$
Mu(x) \leq C u(x), \,\,\,\, \text{a.e.} \,\, x\in \R^d.
$$
If $p<1$ there are no weights so that $M : L^p(u)\to L^{p,\infty}(u)$  is bounded \cite{crs:crs}.

\medskip

Later on, the development of the interpolation theory motivated the study of the boundedness of $M$ on the so called Lorentz spaces.
The (classical) Lorentz space $\llo$ is defined as the class of all functions satisfying  
\begin{align*}
 ||f||_{\Lambda^{p}(w)}:=\left( \int_0^{\infty}  p \, t^{p-1} W(|\{x\in \R^d: |f(x)|>t \}|)dt\right) ^ {1/p}< \infty, 
\end{align*}
where $w$ is a weight in $\R^+$,  $\displaystyle W(t)=\int_0^t w(s)ds$ and $|E|$ denotes the Lebesgue measure of $E$. 
The weak-type Lorentz space $\lloin$ is defined by the following quasi norm
$$
 ||f||_{\Lambda^{p, \infty}(w)}:=\sup_{t>0} t \, W(|\{x\in \R^d: |f(x)|>t \}|)^{\frac 1p}< \infty.
$$
Ari{\~n}o  and Muckenhoupt characterized in  \cite{am:am} the boundedness of $M$ 
on $\llo $. The key idea to study the boundedness of $M$ on these spaces is the existence of $c,C>0$ such that
\begin{equation} \label{lo mismo}
cPf^*(t) \leq (Mf)^*(t) \leq C Pf^*(t).
\end{equation}
In these inequalities $f^*$ is the decreasing rearrangement of $f$, which is defined in $[0,+\infty)$ by 
$$
f^*(t)=\inf\big\{ s>0: |\{ x\in \mathbb R^d: |f(x)|>s\}|\le t \big\},
$$
and $P$ is the Hardy operator defined by
\begin{equation} \label{Hardy}
Pf(t)= \frac{1}{t}\int_0^t f(s) ds, \,\,\,\, t>0,
\end{equation} 
(see \cite{bs:bs} for more details). Consequently, the boundedness of $M$ on $\llo$ is equivalent to the boundedness of 
$P$ on the cone of decreasing functions of $L^p(w)$.
Given $p>0$, the class of weights satisfying
 $$
 M:\llo \to \llo
 $$ 
 is known as $B_p$, and it can be proved \cite{am:am} that $w\in B_p$ if and only if 
\begin{equation*}
r^p \int_r^{\infty}\frac{w(t)}{t^p} \, dt  \le C  \int_0^rw(s) ds, \,\,\,\,\,\, \text{for every}  \,\,\, r>0.
\end{equation*}
Moreover, for every $p>0$, the condition $B_{p,\infty}$ characterizes the boundedness
$$
M:\Lambda^p(w) \longrightarrow  \Lambda^{p, \infty}(w),
$$
where for $p>1$, $B_{p, \infty}=B_p$, and for $p\leq 1$ a weight $w\in B_{p, \infty}$ if and only if 
$$
\frac{W(t)}{t^p} \le C \frac{W(r)}{r^p}, \,\,\,\, \text{for every} \,\,  0<r<t<\infty .
$$
These classes of weights have been well studied in  \cite{am:am, crs:crs, n:n}.

\medskip

 Some analogies between the boundedness properties of $M$ in $L^p(u)$ and in $\llo$ suggested that there might be a unifying theory behind. A natural framework for this unification is provided by the  weighted Lorentz spaces defined by Lorentz in \cite{l1:l1, l2:l2}.  Given $u$,  a weight in $\mathbb R^d$ and given a weight $w$ in $\mathbb R^+$, 
$$
\Lambda^{p}_{u}(w) =\left\{f\in\mathcal M:  ||f||^p_{\Lambda^{p}_{u}(w)}:=\int_0^{\infty}  p  t^{p-1} W(u(\{x\in \R^d: |f(x)|>t \}))dt
< \infty \right\},
$$
where $\mathcal M=\mathcal M(\mathbb R^d)$ is the set of Lebesgue measurable functions on $\mathbb R^d$,  
and the weak-type Lorentz space is defined as follows
$$
\Lambda^{p, \infty}_{u}(w) =\left\{f\in\mathcal M :  ||f||^p_{\Lambda^{p, \infty}_{u}(w)}:=\sup_{t>0} t W^{1/p}(u(\{x\in \R^d: |f(x)|>t \}))< \infty \right\}.
$$
Note that these spaces  include, as particular examples,  the weighted Lebesgue spaces $L^p(u)$,  $L^{p,\infty}(u)$ (when $w=1$) and the Lorentz spaces $\Lambda^p(w)$, $\lloin$ (when $u=1$).

\medskip

In \cite{crs:crs} the strong-type boundedness
\begin{equation} \label{fuerte}
M:\lo \to \lo
\end{equation}
was completely characterized as follows.

\begin{theorem}[\cite{crs:crs}, Theorem 3.3.5]\label{strongmaximal}
For every $0<p<\infty$, 
$$
M:\Lambda^{p}_{u}(w)\longrightarrow \Lambda^{p}_{u}(w)
$$
is bounded if and only if 
 there exists $q\in(0,p)$ such that,  for every finite family of  cubes $(Q_j)_{j=1}^J$,
and every family of measurable sets $(S_j)_{j=1}^{J}$, with $S_j\subset Q_j$, for every $j$,  we have that
\begin{equation}\label{raposo}
\frac{W\left(u\left(\bigcup_{j=1}^J Q_j\right)\right)}{W\left(u\left(\bigcup_{j=1}^J S_j\right)\right)}
       \leq C \max_{1\leq j\leq J} \left(\frac{|Q_j|}{|S_j|}\right)^q.
\end{equation}
for some universal positive constant $C$ depending only on $p$ and the dimension.
\end{theorem}
\noindent It is easy to see that condition \eqref{raposo} recovers $u\in A_p$ if $w=1$, and $w\in B_p$ if $u=1$. Later on, Lerner and P\'erez found  in \cite{lp:lp} other equivalent conditions to the strong boundedness of $M$ in $\lo$ in terms of the so called local maximal operator. 

\medskip

In \cite{crs:crs}, the weak-type boundedness of $M$ was also characterized for $p\leq 1$. In this case, the solution is given by condition \eqref{raposo}, but with the exponent $p$ instead of $q$. However,  the weak-type boundedness 
\begin{equation} \label{debil}
M:\lo \to \loin
\end{equation}
remained open for $p>1$.  The main result in this paper is the following theorem that completely solves this problem.  
\begin{theorem}\label{David y Goliath}
 If    $p>1$,  then
$$
M:\Lambda^{p}_u(w)\longrightarrow \Lambda^{p, \infty}_u(w)
$$
is bounded if and only if \eqref{raposo} holds. In particular, 
$$
M:\Lambda^{p}_u(w)\rightarrow \Lambda^{p, \infty}_u(w)\, \mbox { is bounded}\quad\iff M:\Lambda^{p}_u(w)\rightarrow \Lambda^{p}_u(w)\, \mbox { is bounded}. 
$$\end{theorem}

Finally, we have to mention that, if $d=1$, Theorem \ref{David y Goliath} was proved in~ \cite{aacs:aacs}, and the proof uses
the explicit construction of a function, which together with the weak-type boundedness lead to the geometric condition \eqref{raposo}.  
Even though this paper is inspired on~ \cite{aacs:aacs}, we have to use a different approach, 
since the same method cannot be extended to the multi-dimensional case. 

\subsection*{Notation}
 As usual, we shall use the symbol $A\lesssim B$ to indicate that there exists a universal constant $C$, 
 independent of all important parameters, such that $A\le C B$. Also $A\approx B$ will indicate that 
$A\lesssim B$ and $B\lesssim A$.   
It is known that   the  space $\Lambda^{p}_{u}(w)$  is a  quasi-normed 
space if and only if $w\in \Delta_2$  (see \cite{cgs:cgs}); that is,
$$
W(2r)\lesssim W(r).
$$
This condition will be assumed all over the paper.

\section{Proof of the main result}

\noindent  This section is devoted to the proof of Theorem \ref{David y Goliath}. In some sense, the strategy of the proof combines ideas of 
\cite{n:n}  and~ \cite{aacs:aacs}. We begin with the following two lemmas.

\begin{lemma}\label{f vs chi_E}
Let us assume that 
\begin{equation}\label{wtb}
M:\lo \to \loin
\end{equation}
is bounded. Then, for every $0<\lambda<1$ and every Borel set  $E\subset \R^d$, 
\begin{equation}\label{eq 1}
\|\,\chi_{\{M\chi_E>\la\}}\ M\chi_E\,\|_\lo^p  \lesssim \Big (1+  \log \frac{1}{\lambda}\Big)\ \|\chi_E\|_\lo^p.
\end{equation}
\end{lemma}

\bdem
Fix $0<\lambda<1$. Then

\begin{align*}
\|\,\chi_{\{M\chi_E>\la\}}\ M\chi_E\,\|_\lo^p & 
= \int_0^{\la} p t^{p-1} W(u( \{x: \chi_{\{M\chi_E>\la\}}(x)\ M\chi_E (x)>t   \} )) dt\\
&\quad+  \int_{\la}^1  p t^{p-1} W(u( \{x: \chi_{\{M\chi_E>\la\}}(x)\ M\chi_E(x) >t   \} )) dt \\
&=I+II.
\end{align*}
On the one hand, note that for $t\leq \la$ we have that
$$
\{x: \chi_{\{M\chi_E>\la\}}(x)\ M\chi_E(x) >t   \} = \{x: M\chi_E (x)>\la  \}.
$$
Hence,  by \eqref{wtb}, 
\begin{align*}
I &= \int_0^{\la} p t^{p-1} W\big(u( \{ M\chi_E >\la  \} )\big) dt
 = \la^p  W(u( \{ M\chi_E >\la  \} )) \\
&\leq ||M\chi_E||_{\loin}^p \lesssim ||\chi_E||_{\lo}^p. 
\end{align*}

On the other hand, 
\begin{align*}
II & \leq  \int_{\la}^1  p t^{p-1} W\big(u( \{ M\chi_E >t   \} )\big) dt 
=   p\int_{\la}^1   t^{p} W\big(u( \{ M\chi_E >t   \} )\big) \frac{dt}{t} \\
& \lesssim \int_{\la}^1   ||\chi_E||_\lo^p \frac{dt}{t} = \log \frac{1}{\la} ||\chi_E||_\lo^p, 
\end{align*}
and the result follows.
\edem

\noindent The proof of the following lemma is motived by a result in \cite{lp:lp}.  It provides the extra decay that we shall need to go from the weak-type to the strong-type boundedness.

\begin{lemma}\label{LP}
For any $0<\la<1$ and any Borel subset $E\subset \R^d$, it holds that
\begin{equation} \label{eq 2}
\chi_{ \{M\chi_E>\la\}}(x)\lesssim \frac{1}{\la(1-\log\la)} M\big(\chi_{ \{M\chi_E>\la\}}\,M\chi_E\big)(x)\quad (x\in\R^d).
\end{equation}
\end{lemma}

\bdem
Fix a Borel set $E\subset \R^d$,  $\lambda \in (0,1)$ and $x\in\R^d$ such that  $M\chi_E(x)>\la$. Then there exists a cube $Q$ so that $x\in Q$ and 
$$
\la<\frac{|E\cap Q|}{|Q|}.
$$
Since  the function 
$
\phi(x)=x\left(1+\log \frac{1}{x}\right)
$
  is     increasing   in $(0,1)$, we have that
\begin{align*}
\la \Big(1+\log \frac{1}{\la}\Big)&=\phi(\la) \leq \phi \left( \frac{|E\cap Q|}{|Q|} \right)
=\frac{1}{|Q|} \int_0^{|Q|} \frac{\min (t, |E\cap Q|)}{t}dt\\
& = \frac{1}{|Q|}\int_0^{|Q|} P(\chi_{E\cap Q})^*(t) dt,
\end{align*}
where $P$ denotes the Hardy operator defined by \eqref{Hardy}. 
Hence, by \eqref{lo mismo},  we obtain
\begin{align*}
\la \Big(1+\log \frac{1}{\la}\Big)&\approx \frac{1}{|Q|}\int_0^{|Q|} (M\chi_{E\cap Q})^*(t) dt \\
&\leq  \frac{1}{|Q|}\int_0^{|Q|}\big( \chi_{3Q }M\chi_{E\cap Q}\big)^*(t) dt 
+ \frac{1}{|Q|}\int_0^{|Q|} \big(\chi_{(3Q)^c }M\chi_{E\cap Q}\big)^*(t) dt\\
 &\leq  \frac{1}{|Q|}\int_{3Q} M\chi_{E\cap Q} (y) dy +  \frac{1}{|Q|}\int_0^{|Q|} \big(\chi_{(3Q)^c }M\chi_{E\cap Q}\big)^*(t) dt.
\end{align*}
Now, the standard estimate
$$
\chi_{(3Q)^c}(z)M\chi_{E\cap Q}(z) \lesssim \inf_{y \in Q} M\chi_{E\cap Q} (y)\le \inf_{y \in Q} M\chi_{E} (y), \quad z\in\R^d
$$
implies that
\begin{align*}
\la \Big(1+\log \frac{1}{\la}\Big)&\lesssim
\frac{1}{|Q|}\int_{3Q} M\chi_{E} (y) dy +  \frac{1}{|Q|}\int_Q M\chi_{E}(y) dy\\
&\lesssim M( M\chi_E)(x) \le M( \chi_{\{M\chi_E>\la\}} M\chi_E)(x) + M( \chi_{\{M\chi_E\le \la\}} M\chi_E)(x)\\
&\le M( \chi_{\{M\chi_E>\la\}} M\chi_E)(x)  + \lambda. 
\end{align*}

Finally, since $\{M\chi_E>\la\}$ is an open set, we obviously have  that 
$$\lambda \le M( \chi_{\{M\chi_E>\la\}} M\chi_E)(x) $$ 
and hence  the result follows. \edem

\noindent Equivalently, we can write the inequality \eqref{eq 2} as an inclusion of level sets in the following way.

\begin{corollary}\label{inclusion}
There  exists $c>0$ such that, for every Borel subset  $E\subset \R^d$  and every $0<\lambda<1$, 
$$
\{M\chi_E>\lambda\}\subseteq \{M( \,\chi_{\{M\chi_E>\la\}}\ M\chi_E \,)>c\,\la(1-\log\la)\}.
$$
\end{corollary}

\noindent Now, in order to proceed to the proof of our main theorem, we need to recall the following result proved in \cite{crs:crs} (see Theorems 3.3.3 and 3.3.5). 

\begin{proposition} \label{JA} If there exists $0<r<\infty$ such that
\begin{align*}
\int_0^1 \lambda^{r-1} W^{r/p}\big(u(\{ M\chi_E>\lambda\})\big) d\lambda\lesssim  ||\chi_E||^p_{\lo}, 
\end{align*}
then \eqref{raposo}  holds. 
\end{proposition}

\bdem[Proof of Theorem \ref{David y Goliath}]
Let $0<\lambda<1$ and $f=\chi_{\{M\chi_E>\la\}}\ M\chi_E$. By Corollary \ref{inclusion},  we have that
$$
W\big(u(\{\,M\chi_E>\la\})\big)\le  W\big(u(\{\,Mf>c \la (1-\log\la)\})\big), 
$$
and using the weak-type boundedness of $M$,  it holds that
$$
W\big(u(\{Mf>c \la (1-\log\la)\})\big)\lesssim \frac{1}{\la^p(1-\log\la)^{p}}\|f\|_{\lo}^p.
$$
By \eqref{eq 1}  we obtain that 
$$
W\big(u(\{\,M\chi_E>\la\})\big)\lesssim \frac{1}{\la^p(1-\log\la)^{p-1}}\|\chi_E\|_{\lo}^p, 
$$
and hence,   if we take $r>0$ such that  $p/(p-1) <r<\infty$,  we have that 
\begin{align*}
 \int_0^1 \lambda^{r-1} W^{r/p}\big(u(\{ M\chi_E>\lambda\})\big) d\lambda \lesssim ||\chi_E||^p_{\lo}
\end{align*} 
and the result follows by Proposition \ref{JA}.
\edem

\bigskip

\proof[Acknowledgements]
We would like to thank Prof.\  Javier Soria for the helpful discussions related with the subject of this paper, and also for his useful comments that let us improve the presentation of this paper. The first and second authors would also like to thank the University of Barcelona (UB) and the Institute of Mathematics of the University of Barcelona (IMUB) for providing us all the facilities and hosting us during the research stay that led to this collaboration.

\bigskip



\end{document}